\documentclass{ifacconf}

\usepackage{graphicx}      
\usepackage{natbib}        

\usepackage{amsmath}
\usepackage{amssymb}
\usepackage{bbold}

\usepackage{color}
\usepackage{comment}
\begin{document}
\begin{frontmatter}

\title{Design of an optimal positive state observer for~a~diffusion~system~with~point~observation~\thanksref{footnoteinfo}} 


\thanks[footnoteinfo]{This research was conducted with the financial support of F.R.S-FNRS. Violaine Piengeon is a FNRS Research Fellow under the grant FC 58057.}

\author[UNamur]{V. Piengeon} 
\author[UNamur]{and J.J. Winkin}

\address[UNamur]{University of Namur, Department of Mathematics \& naXys Institute, Rue Graf\'e 2, B-5000 Namur, Belgium\\ (e-mails: violaine.piengeon@unamur.be, joseph.winkin@unamur.be).}

\begin{abstract}                
An optimal state observer design problem for a diffusion system with point observation is studied. As the considered system is positive, the analysis takes into account this additional condition. The solution of the optimization problem is derived by means of a discretization approach, which consists in considering the analogous problem adapted for a spatially discretized version of the pure diffusion system. 
\end{abstract}


\end{frontmatter}

\section{Introduction}
In many contexts such as feedback control, it is necessary that the state of a dynamical system be sufficiently well-known. However, in practice, it is usually impossible to measure the full state. This limitation is overcome by using state estimation, which aims to reconstruct the trajectories of the real system by means of available measurements, via appropriate sensors. This gives rise to the design of observers. For linear distributed parameter systems, \cite{Curtain_Zwart2020} and \cite{Tello2021} provide an overview of existing observers.

Here, we focus on the design of observers for a specific class of distributed parameter systems, namely positive systems, which have the particularity to leave the nonnegative cone of the state space invariant; see e.g. \cite{Smith1995}, \cite{Batkai2017}. Positive systems have been extensively studied because the positivity property is relevant in many fields such as (bio)chemical engineering, biomedical sciences, or epidemiology; see e.g.  \cite{Laabissi2001}, \cite{SONVEAUX2023}, \cite{Dochain2025}. Indeed, physical quantities such as temperature (K), concentration or population density are all nonnegative quantities that can be represented by state variables. When one wants to estimate the state of a positive system, it is necessary to design state observers that produce estimated states maintaining this natural positivity property, otherwise the estimated state would not have any physical meaning. The literature on observers for positive linear finite-dimensional systems is quite extensive, see e.g. \cite{Rami2011} and references therein. 

In this work, we consider the state estimation problem for a specific positive system, which is the (pure) diffusion system described by the partial differential equation (PDE)
\begin{equation}\label{eq:pure_diffusion}
    \frac{\partial x}{\partial t}(z,t) = D_a\frac{\partial^2 x}{\partial {z}^2}(z,t), \quad z\in [0,L]\, ,\: t\geq 0 ,
\end{equation}
with Neumann boundary conditions 
\begin{equation}\label{eq:cdb}
 \dfrac{\partial x}{\partial z}(0,t)  =  0, \quad 
\dfrac{\partial x}{\partial z}(L,t)  =  0 ,
\end{equation}
where~$D_a$ and $L>0$ are positive constants and $x(.,t)$ is the state at time $t$. An initial condition is of the form ~$x(z,0) = x_0(z), \;\; z\in [0,L]$. The measured output $y$ that is considered consists in a boundary observation given by
\begin{equation}\label{eq:output}
    y(t) = x(L,t), \quad t\geq 0.
\end{equation}
Note that the PDE~\eqref{eq:pure_diffusion} can be seen as the well-known heat equation on the spatial domain $[0,L]$, with diffusion coefficient $D_a$; see e.g. \cite[Subsection 6.1.3]{VandeWouwer2014} or \cite[Section 1]{Crank1975}. This diffusion system belongs to the class of positive systems: indeed, the associated operator defined by $A_{0} = D_a({\rm  d}^2/{\rm  d}  z^2)$ with domain
\begin{equation*}\label{eq:operator_A0}
    D(A_0) = \left\{x\in \mathrm{H}^2(0,L) : \dfrac{{\rm  d}  x}{{\rm  d}  z}(0) =  0 = \dfrac{{\rm  d}  x}{{\rm  d}  z}(L)  \right\},
\end{equation*}
generates a positive $C_0$-semigroup on the state space $\mathrm{L}^2(0,L)$; see \cite{Laabissi2001}. Here, $\mathrm{H}^2(0,L)$ denotes the Sobolev space of functions belonging to the Lebesgue space $\mathrm{L}^2(0,L)$ whose first and second weak derivatives are also in $\mathrm{L}^2(0,L)$.

The present work aims at providing the solution of an optimal (according to a well-chosen optimality criterion) state observer design problem for the diffusion system \eqref{eq:pure_diffusion}-\eqref{eq:output}. This problem (see Section \ref{section:Optimal state estimation problem}) may be interpreted as a trade-off between the advantages and an attenuation of the disadvantages of a high-gain type observer. Here, the positivity property of the observer is a constraint of the optimal state estimation problem so that the optimal solution is guaranteed to satisfy this condition. Our approach is based on the following scheme: discretization of the system with respect to the spatial variable (using a finite difference scheme), solution for the discretized model (Theorem \ref{theorem:design_k1_opti}), taking the limit as the discretization step tends to zero (Lemma \ref{lemma:conv_sol}) and validation for the nominal system (Theorem \ref{theorem:design_gamma_opti}). This methodology and the revealed links between the nominal problem and the discretized one constitute a key feature of our work. Moreover, the analysis of the discretized system is interesting in itself because it provides a procedure to find an analytical and computable expression of the related optimal solution (see Section \ref{section:Discretization approach}).

\section{Positive State Observers}\label{section:estimation}
In order to get an estimate of the exact state, we consider the state observer for the diffusion system \eqref{eq:pure_diffusion}-\eqref{eq:output}, which is based on the predictor-corrector scheme given by 
\begin{equation} \label{eq:pure_diffusion_observer}
    \frac{\partial \tilde{x}}{\partial t}(z,t) = D_a\frac{\partial^2 \tilde{x}}{\partial z^2}(z,t), \quad z \in [0,L], \: t \geq 0,
\end{equation}
with the boundary conditions
\begin{equation} \label{eq:observer_BC}
        \frac{\partial \tilde{x}}{\partial z}(0,t) = 0, \quad
        \frac{\partial \tilde{x}}{\partial z}(L,t) - \gamma(\tilde{y}(t) - y(t)) = 0,
\end{equation}
where the measured output $y$ is given by \eqref{eq:output}, and where
\begin{equation}\label{eq:estimated_output}
    \tilde{y}(t):=\tilde{x}(L,t),\quad t\geq 0,
\end{equation}
denotes the estimated output, and $\gamma \in \mathbb{R}$ is the observer gain, to be designed. The dynamics of the estimation error $e:=\tilde{x}-x$ is described by the following PDE obtained by subtracting {\eqref{eq:pure_diffusion}-\eqref{eq:cdb} from \eqref{eq:pure_diffusion_observer}-\eqref{eq:observer_BC}:
\begin{equation} \label{eq:error_dynamics}
    \frac{\partial e}{\partial t}(z,t) = D_a\frac{\partial^2 e}{\partial z^2}(z,t), \quad z \in [0,L], \: t \geq 0,
\end{equation}
with the boundary conditions
\begin{equation} \label{eq:error_BC}
        \frac{\partial e}{\partial z}(0,t) = 0, \quad
        \frac{\partial e}{\partial z}(L,t) - \gamma e(L,t) = 0.
\end{equation}
By considering the invertible map $S : \mathrm{H}^2(0,L)\to \mathrm{H}^2(0,L)$ given by $(Sx)(z) = x(L-z)=:e(z)$, we observe that the estimation error dynamics \eqref{eq:error_dynamics}-\eqref{eq:error_BC}
is equivalent to the closed-loop system analyzed in \cite{Dehaye2016},
where $\gamma$ replaces the control gain $-\kappa$. Consequently, the estimation error dynamics~\eqref{eq:error_dynamics}-\eqref{eq:error_BC} is exponentially stable, or equivalently the operator 
$A_{\gamma} = D_a({\rm  d}^2/{\rm  d}  z^2)$, with domain
\begin{equation*}\label{eq:operator_Ag}
   D(A_{\gamma}) = \left\{e\in \mathrm{H}^2(0,L) : \dfrac{{\rm  d}  e}{{\rm  d}  z}(0) = 0 = \dfrac{{\rm  d}  e}{{\rm  d}  z}(L) - \gamma e(L) \right\},
\end{equation*}
is the infinitesimal generator of an exponentially stable $C_0$-semigroup on $\mathrm{L}^2(0,L)$ if, and only if, $\gamma < 0$; see \cite[Theorem 4]{Dehaye2016}. In what follows, the $C_0$-semigroup generated by $A_{\gamma}$ on $\mathrm{L}^2(0,L)$ and the associated growth bound are denoted by $(T_{A_{\gamma}}(t))_{t\geq0}$ and $\omega_0^{\gamma}$, respectively. 

Beyond stability of the estimation error dynamics, it is essential to ensure the nonnegativity invariance property of the estimated state. In this context, the following result holds. 
\begin{prop}\label{prop:positivity_A0gamma}
    For all $\gamma< 0$, the operator $\mathcal{A}_{\gamma}$ defined by
    \begin{equation*}
        \mathcal{A}_{\gamma} = \begin{pmatrix}
            D_a({\rm  d}^2/{\rm  d}  z^2) & 0\\
            0 & D_a({\rm  d}^2/{\rm  d}  z^2)
        \end{pmatrix},
    \end{equation*}
    with domain
    \begin{equation*}\label{eq:operator_A0gamma}
    \begin{array}{rcl}
            D(\mathcal{A}_{\gamma}) &=& \left\{(x,\tilde{x})\in \mathrm{H}^2(0,L) \times \mathrm{H}^2(0,L):\right.\\ &&\dfrac{{\rm  d}}{{\rm  d}  z}\begin{pmatrix}x\\\tilde{x}\end{pmatrix}(0) = 0,\\ &&\left.\dfrac{{\rm  d}}{{\rm  d}  z}\begin{pmatrix}x\\\tilde{x}\end{pmatrix}(L) = \begin{pmatrix}
            0&0\\-\gamma&\gamma\end{pmatrix}\begin{pmatrix}x\\\tilde{x}\end{pmatrix}(L) \right\},  
    \end{array}
    \end{equation*}
    generates a positive $C_0$-semigroup~$(T_{\mathcal{A}_{\gamma}}(t))_{t\geq 0}$ on~$\mathrm{L}^2(0,L) \times \mathrm{L}^2(0,L)$. 
\end{prop}
A consequence of Proposition \ref{prop:positivity_A0gamma} is that, by combining the state dynamics \eqref{eq:pure_diffusion}-\eqref{eq:cdb} and the state observer dynamics \eqref{eq:pure_diffusion_observer}-\eqref{eq:observer_BC} with $\gamma<0$, the resulting system is positive. 
Thus, for every nonnegative initial state $x_0$ and every nonnegative estimated initial state $\tilde{x}_0$ such that $(x_0,\tilde{x}_0)\in D(\mathcal{A}_{\gamma})$, the associated estimated state trajectory $\tilde{x}$ is nonnegative when $\gamma<0$. 
Therefore, the proposed state observer \eqref{eq:pure_diffusion_observer}-\eqref{eq:observer_BC} is well-suited to positive state estimation if, and only if, $\gamma < 0$.

\section{Optimal state estimation problem}\label{section:Optimal state estimation problem}
Since there exists an infinite number of state observers of the form \eqref{eq:pure_diffusion_observer}-\eqref{eq:observer_BC} that meet the positive state estimation criteria, the goal is to find the best one according to a well-chosen optimality criterion. More precisely, let us consider, for any $\alpha>0$, the optimal state estimation problem for the pure diffusion system, denoted $\mathcal{S}_{\alpha}$, given by
\begin{subequations}\label{eq:pbm_opti}
    \begin{align}
        \displaystyle \min_{\gamma \,\in \,\mathbb{R}}\quad &\gamma^2, \label{eq:pbm_opti_obj}\\
        \text{s.t.}\quad & (T_{\mathcal{A}_{\gamma}}(t))_{t\geq0} \text{ is positive,} \label{eq:pbm_opti_const_positivity}\\
        \quad & (T_{A_{\gamma}}(t))_{t\geq0} \text{ is } \beta\text{-exponentially stable},\notag\\ 
        \quad & \,\,\forall \beta\in(-\alpha,0].\label{eq:pbm_opti_const_stability}
    \end{align}
\end{subequations}
\begin{rem}
$\alpha$) Condition~\eqref{eq:pbm_opti_const_stability} is a stability constraint which is stronger than exponential stability: it is equivalent to $\omega_0^{\gamma}\leq-\alpha$; see \cite[Section 4.1]{Curtain_Zwart2020}. The latter non-strict inequality ensures the feasibility of the optimal state estimation problem $\mathcal{S}_{\alpha}$. 
 In condition \eqref{eq:pbm_opti_const_stability}, the parameter $\alpha$ corresponds to the required minimum stability margin for $(T_{A_{\gamma}}(t))_{t\geq0}$. \\
$\beta$) A state observer of the form~\eqref{eq:pure_diffusion_observer}-\eqref{eq:observer_BC} where $\gamma$ satisfies~\eqref{eq:pbm_opti_const_positivity} and~\eqref{eq:pbm_opti_const_stability} may be viewed as a constrained high gain (static) state observer. The high gain interpretation comes from the fact that the output error injection term $\gamma(\tilde{y}(t) - y(t))$
is parametrized by a single scalar $\gamma$ and can become arbitrarily large 
as $|\gamma|\rightarrow \infty$. This kind of observer is widely used in the literature: see \cite{Khalil2008} and references therein or \cite{Kitsos2021}. In these references, it is mentioned that a key feature of high gain type observers is a fast convergence of the estimation error. Here, this effect is sought with condition \eqref{eq:pbm_opti_const_stability} and the choice of an $\alpha$ sufficiently large (in the following, we will see that $\alpha$ cannot be driven too high). However, high gain type observers may be subject to the peaking phenomenon; see \cite{Khalil2008}. Numerical simulations that we performed  for the state observer \eqref{eq:pure_diffusion_observer}-\eqref{eq:observer_BC} reveal also this phenomenon. In this context, the optimality criterion~\eqref{eq:pbm_opti_obj} aims to prevent this effect.

\end{rem}

In the light of those comments, the $\mathcal{S}_{\alpha}$ problem may be viewed as a trade-off between the advantage offered by a sufficiently high gain state observer (in terms of decay rate of the estimation error dynamics) and the attenuation of the peaking phenomenon, while ensuring nonnegativity of the estimated state trajectories.

\section{Discretization approach}\label{section:Discretization approach}
With a view to solving the $\mathcal{S}_{\alpha}$ problem, for any $\alpha>0$, both theoretically and numerically, the analysis is based on an optimal state estimation problem associated to spatially discretized versions of the considered dynamics. 

Before reporting the studied optimization problem, we introduce the approximations under consideration. We consider $n$ discretization points $z_i$, $i=1,\dots,n$, with $z_1=0$ and $z_n=L$. The discretization step, denoted by $\Delta z$, is given by~$\Delta z = L/(n-1)$.
Choosing the same finite difference schemes and the same notations than \cite[Section 3.2]{Dehaye2016}, we obtain the discretized pure diffusion system 
\begin{equation}\label{eq:pure_diffusion_discr}
    \dot{x}^{(n)}(t) = A^{(n)}x^{(n)}(t)\,,\quad t\geq 0\,,
\end{equation}
and the discretized state observer
\begin{equation}\label{eq:observer_discr}
    \dot{\tilde{x}}^{(n)}(t) = A^{(n)}\tilde{x}^{(n)}(t) +l^{(n)}(\tilde{y}(t)-y(t))\,,\quad t\geq 0\,,
\end{equation}
where 
\begin{equation*}\label{eq:matrice_An}
    A^{(n)} = 
    \begin{pmatrix}
        -p_2 & p_2 & 0 & \cdots  & 0\\
        p_2 & -2p_2 & p_2 &  & \vdots\\
        0 & \ddots & \ddots & \ddots  & 0\\
        \vdots & & p_2 & -2p_2 & p_2\\
        0  & \cdots & 0 & p_2 & -p_2
    \end{pmatrix}\in {\mathbb R}^{n\times n},
\end{equation*}
is a tridiagonal matrix, 
$$l^{(n)} = \gamma p_1 \mathbb{e}_n^{(n)} \in {\mathbb R}^n ,$$
$\mathbb{e}_n^{(n)}$ is the last vector of the canonical basis of ${\mathbb R}^n$, namely $ \mathbb{e}_n^{(n)} := \begin{pmatrix} \; \delta_{ni} \;  \end{pmatrix}_{i=1, \cdots,n}^\top$, where $\delta_{kl}$ is the usual Kronecker symbol, and
\[
p_1 = \dfrac{D_a}{\Delta z},\quad p_2 = \dfrac{D_a}{\Delta z^2}.
\]
The state trajectories of~\eqref{eq:pure_diffusion_discr} and \eqref{eq:observer_discr}
are defined by
\begin{equation*}\label{eq:vecteur_xn}
    x^{(n)}(\cdot) :=
    \begin{pmatrix}
        x(z_1,\cdot) & \cdots & x(z_n,\cdot)   
    \end{pmatrix}^{\top} \in {\mathbb R}^n.
\end{equation*}
and 
\begin{equation*}\label{eq:vecteur_xtilde_n}
    \tilde{x}^{(n)}(\cdot) :=
    \begin{pmatrix}
        \tilde{x}(z_1,\cdot) & \cdots & \tilde{x}(z_n,\cdot)   
    \end{pmatrix}^{\top} \in {\mathbb R}^n,
\end{equation*}
respectively. Moreover, the measured and estimated outputs can be written as 
\[
    y(t) = C^{(n)}x^{(n)}(t) \text{ and } \tilde{y}(t) = C^{(n)}\tilde{x}^{(n)}(t),\quad t\geq 0,
\]
where $$C^{(n)} ={\mathbb{e}_n^{(n)}}^{\top}\in {\mathbb R}^{1\times n}.$$
The dynamics of the associated estimation error, defined by $e^{(n)}:=\tilde{x}^{(n)}-x^{(n)}$, is given by 
\begin{equation}\label{eq:error_discr}
    \dot{e}^{(n)}(t) = A^{(n)}_{\gamma}e^{(n)}(t)\,,\quad t\geq 0\,,
\end{equation}
where 
\begin{equation*}
    A^{(n)}_{\gamma} := A^{(n)}+l^{(n)}C^{(n)},
\end{equation*}
and corresponds to a spatially discretized version of the nominal error dynamics. 

One of the advantages of considering this discretization method is that the obtained discretized dynamics shares common features with the nominal system. Indeed, the discretized diffusion system \eqref{eq:pure_diffusion_discr} is a positive LTI system, that is, $A^{(n)}$ is a Metzler matrix and $C^{(n)}$ is a nonnegative matrix. Observe also that, in addition to being tridiagonal, the matrix $A^{(n)}$ is symmetric, and so all its eigenvalues are real; in addition, zero is its unique unstable eigenvalue. For $n$ sufficiently large, the estimation error dynamics is exponentially stable, that is the matrix $A^{(n)}_{\gamma}$ is exponentially stable, if, and only if,~$\gamma< 0$; see \cite[Theorem 4]{Dehaye2016}. Finally, if $\gamma< 0$, then it is easy to see that the extended system \eqref{eq:pure_diffusion_discr}-\eqref{eq:observer_discr} is positive, or equivalently the following matrix is Metzler:
\begin{equation}
    \mathcal{A}_{\gamma}^{(n)} := \begin{pmatrix}
        A^{(n)} & 0\\
        -l^{(n)}C^{(n)} & A^{(n)}_{\gamma}
    \end{pmatrix}.
\end{equation}

Let us now consider, for any $\alpha>0$, the optimal state estimation problem for the discretized pure diffusion system, denoted  $\mathcal{S}_{\alpha}^{(n)}$, given by
\begin{subequations}\label{eq:pbm_opti_discr}
    \begin{align}
        \displaystyle \min_{\gamma \,\in \,\mathbb{R}}\quad &{\gamma}^2, \label{eq:pbm_opti_discr_obj}\\
        \text{s.t.}\quad & \mathcal{A}_{\gamma}^{(n)} \text{ is Metzler,} \label{eq:pbm_opti_discr_const_positivity}\\
        \quad & A^{(n)}_{\gamma} \text{ is } \beta\text{-exponentially stable}, \forall \beta\in(-\alpha,0]\,,\label{eq:pbm_opti_discr_const_stability}
    \end{align}
\end{subequations}
For all~$\gamma\in\mathbb{R}$, since $A^{(n)}_{\gamma}$ is symmetric, all its eigenvalues are real. Therefore, the  constraint~\eqref{eq:pbm_opti_discr_const_stability} is equivalent to the condition $\max \{\lambda : \lambda\in\sigma(A^{(n)}_{\gamma})\} \leq-\alpha$.\\
As for the $\mathcal{S}_{\alpha}$ problem, the feasible set of the $\mathcal{S}_{\alpha}^{(n)}$ problem is the set of all values of $\gamma$ that yield a positive state observer and an exponentially stable estimation error dynamics with an associated stability margin greater or equal to $\alpha$. Moreover, if the $\mathcal{S}_{\alpha}^{(n)}$ problem is feasible, then the pair $(C^{(n)},A^{(n)})$ is $\beta$-exponentially detectable for all $\beta\in(-\alpha,0]$; see \cite[Section 8.1]{Curtain_Zwart2020}.

The next lemma gives a simple necessary and sufficient condition on $\alpha$ for the $\mathcal{S}_{\alpha}^{(n)}$ problem to be feasible. It contains also an algebraic scalar inequality, which is equivalent to the combination of both constraints~\eqref{eq:pbm_opti_discr_const_positivity} and~\eqref{eq:pbm_opti_discr_const_stability}. 
\begin{lem}\label{lemma:feasability_practicality}
    For all ~$\alpha>0$, the $\mathcal{S}_{\alpha}^{(n)}$ problem~\eqref{eq:pbm_opti_discr} for the discretized  diffusion system \eqref{eq:pure_diffusion_discr} is feasible, i.e. there exists $\gamma$ such that \eqref{eq:pbm_opti_discr_const_positivity} and~\eqref{eq:pbm_opti_discr_const_stability} hold, if, and only if,~$\alpha < -\alpha_f^{(n)}$, where 
    \begin{equation*}
        \alpha_f^{(n)} := \max \{\lambda \in \mathbb{R}: |A^{(n)}-\lambda I|_{n-1} = 0 \},
    \end{equation*}
    and the notation $|\cdot|_j$ is used for the leading principal minor of order $j$, for any integer $j$ such that $1\leq j \leq n$. Moreover, for all $\gamma\in\mathbb{R}$, conditions~\eqref{eq:pbm_opti_discr_const_positivity} and~\eqref{eq:pbm_opti_discr_const_stability} hold if, and only if,
    \begin{equation}\label{eq:designk1}
        \gamma\leq\frac{p_2-\alpha}{p_1}+\frac{p_2^2}{p_1}\frac{|A^{(n)}+\alpha I|_{n-2}}{|A^{(n)}+\alpha I|_{n-1}}.
    \end{equation}
\end{lem}
\vspace{1mm}
\begin{rem}
The feasibility condition on $\alpha$ in Lemma \ref{lemma:feasability_practicality}
means that there are no values of $\gamma$ that give a stability margin associated to $A^{(n)}_{\gamma}$ that is greater or equal to $-\alpha_f^{(n)}$: it is not possible to ensure any (i.e. too large) stability margin for the estimation error dynamics \eqref{eq:error_discr} while ensuring positivity of the discretized state observer \eqref{eq:observer_discr}. The condition $\alpha < -\alpha_f^{(n)}$ is also a sufficient condition for the $\beta$-exponential detectability of $(C^{(n)},A^{(n)})$ for all $\beta\in(-\alpha,0]$. 
\end{rem}
For any $n$, the following result provides an analytical and explicit solution of the $\mathcal{S}_{\alpha}^{(n)}$ problem when it is feasible. 
\begin{thm}\label{theorem:design_k1_opti}
    Fix~$\alpha\in\mathbb{R}$ such that $\,\alpha_f^{(n)} <  -\alpha < 0$. 
The $\mathcal{S}_{\alpha}^{(n)}$ problem~\eqref{eq:pbm_opti_discr} for the discretized diffusion system \eqref{eq:pure_diffusion_discr} has a unique solution~$\gamma^{(n)}_*$, which is given by
    \begin{equation}\label{eq:solution_discr}
        \gamma^{(n)}_* =\frac{\sqrt{D_a}}{D_a} \cdot
        \frac{\displaystyle\sum_{l=0}^{n-1} \binom{n+l}{2l+1}\frac{(-1)^{l+1}}{p_2^{l+\frac{1}{2}}}\alpha^{l+1}}{\displaystyle\sum_{l=0}^{n-1} \binom{n+l-1}{2l}\frac{(-1)^{l}} {p_2^{l}}\alpha^{l}}.
    \end{equation}
\end{thm}
The first step in the proof of Theorem \ref{theorem:design_k1_opti} consists in showing that the right-hand side of \eqref{eq:designk1} is negative. Thus, in view of Lemma \ref{lemma:feasability_practicality}, the $\mathcal{S}_{\alpha}^{(n)}$ problem \eqref{eq:pbm_opti_discr} has a unique solution~$\gamma^{(n)}_*$, which is first given by 
\begin{equation}\label{eq:solution_discr_exp1}
    \gamma^{(n)}_* = \frac{p_2-\alpha}{p_1}+\frac{p_2^2}{p_1}\frac{|A^{(n)}+\alpha I|_{n-2}}{|A^{(n)}+\alpha I|_{n-1}}.
\end{equation}
However, this expression of the solution is not explicit. So, we derived \eqref{eq:solution_discr} from \eqref{eq:solution_discr_exp1} in order to get a computable expression that is easier to use in the convergence analysis that is reported hereafter, \cite{Piengeon2026}.  

\section{Convergence and validation}
To infer a solution of the optimal state estimation problem associated to the pure diffusion system, the limit of the solution of the discretized problem can be computed.
\begin{lem}\label{lemma:conv_sol}
    Fix~$\alpha\in\mathbb{R}$ such that $0< \alpha<\pi^2 D_a/(4L^2)$. The solution $\gamma^{(n)}_*$ of the $\mathcal{S}_{\alpha}^{(n)}$ problem~\eqref{eq:pbm_opti_discr} converges asymptotically to $\gamma_*$ given by~
    \begin{equation}\label{eq:solution}
        \gamma_* :=  -\sqrt{\frac{\alpha} {D_a}}\tan{\left(L\cdot\sqrt{\frac{\alpha}{D_a}}\right)},
    \end{equation}
    i.e.  
    \begin{equation*}\label{eq:convergence}
     \gamma^{(n)}_*\mathcal{X}_{\left(0,-\alpha_f^{(n)}\right)}(\alpha) \rightarrow   \gamma_*\quad \text{as}\quad n\rightarrow \infty\,.
    \end{equation*}
\end{lem}
The results obtained for the discretized model in Theorem~\ref{theorem:design_k1_opti} suggest that the limit computed in Lemma~\ref{lemma:conv_sol} might be the solution of the optimal state estimation problem associated to the nominal pure diffusion system. This conjecture may actually be validated. 
\begin{thm}\label{theorem:design_gamma_opti}
    Fix~$\alpha>0$. The $\mathcal{S}_{\alpha}$ problem~\eqref{eq:pbm_opti} for the diffusion system \eqref{eq:pure_diffusion}-\eqref{eq:cdb} is feasible if, and only if, 
    $\alpha <\pi^2 D_a/(4L^2)$. In this case, the $\mathcal{S}_{\alpha}$ problem has a unique solution $\gamma_*$ which is given by \eqref{eq:solution}.
\end{thm}
It follows from Theorem~\ref{theorem:design_gamma_opti} that it is not possible to design a positive state observer of the form~\eqref{eq:pure_diffusion_observer}-\eqref{eq:observer_BC} for the diffusion system \eqref{eq:pure_diffusion}-\eqref{eq:cdb} with a stability margin associated to the estimation error dynamics which is greater or equal to $\pi^2 D_a/(4L^2)$. Following a remark done for the discretized case, the feasibility condition of the $\mathcal{S}_{\alpha}$ problem can also be interpreted as a sufficient condition for the $\beta$-exponential detectability of the considered diffusion system with the boundary observation \eqref{eq:output}. Furthermore, if~$0<\alpha <\pi^2 D_a/(4L^2)$, then $|\gamma_*|$ is the smallest value of $|\gamma|$ that yields a positive state observer associated to an estimation error dynamics with a stability margin greater or equal to~$\alpha$.

\section{Concluding remarks}
A class of positive state observers for the considered pure diffusion system with point observation is identified. Among these, we have found the best one according to a trade-off between a sufficiently high stability margin for the estimation error dynamics and the reduction of the so-called peaking phenomenon. We have also derived a feasibility condition for this optimization problem and observed that it may be interpreted as a sufficient condition for $\beta$-detectability, for all $\beta\in(-\alpha,0]$, of the considered system. The discovery of the solution of the nominal problem is presented via a detailed study of a spatially discretized version of the problem. This methodology allows to exhibit the links between the nominal problem and its discretized version.

We have also addressed the dual problem of optimal positively stabilizing state feedback design for the pure diffusion system with point control, as well as the design of positively stabilizing observer-based state feedbacks,  
\cite{Piengeon2026}. An interesting open question is the adaptation of the discretization approach to exhibit the solution of similar optimal state observer design problems for other positive distributed systems, such as convection-diffusion-reaction systems with more general boundary conditions.



\bibliography{ifacconf}             
                                                   







\end{document}